\documentclass[12pt]{article}

\usepackage{amsmath}
\usepackage{amssymb}

\usepackage{times}
\usepackage{bm}
\usepackage{natbib}

\usepackage{paralist} 
\usepackage{mathrsfs}
\usepackage{multirow}
\usepackage{graphicx} 
\usepackage{subfig}
\usepackage{enumitem}
\usepackage{array,booktabs}
\usepackage{float}

\usepackage{authblk}
\usepackage{amsthm}
\newtheorem{theorem}{Theorem}
\newtheorem{lemma}{Lemma}
\newtheorem{corollary}{Corollary}
\newtheorem{definition}{Definition}
\newtheorem{remark}{Remark}

\usepackage[toc,page,title,titletoc,header]{appendix}

\newcommand{\essinf}{{\mathrm{ess}\inf}}

\newcommand{\Proj}{{\mathbf P}}
\newcommand{\rd}{\,\mathrm{d}}
\newcommand{\EP}{\,\mathbb{P}}
\newcommand{\EE}{\,\mathbb{E}}

\newcommand{\bsbX}{{\boldsymbol{X}}}
\newcommand{\bsbx}{{\boldsymbol{x}}}
\newcommand{\bsby}{{\boldsymbol{y}}}
\newcommand{\bsbb}{{\boldsymbol{\beta}}}
\newcommand{\bsbg}{{\boldsymbol{\gamma}}}

\newcommand{\bsbI}{{\boldsymbol{I}}}

\newcommand{\bsbDelta}{{\boldsymbol{\Delta}}}

\newcommand{\bsbxi}{{\boldsymbol{\xi}}}

\newcommand{\bsba}{{\boldsymbol{\alpha}}}
\newcommand{\bsbA}{{\boldsymbol{A}}}

\newcommand{\bsbeps}{{\boldsymbol{\epsilon}}}

\newcommand{\bsbh}{{\boldsymbol{h}}}
\newcommand{\bsbd}{{\boldsymbol{d}}}

\DeclareMathOperator{\vect}{\mbox{vec}\,}

\usepackage[normalem]{ulem}

\begin{document}

\title{On the Finite-Sample Analysis of   $\Theta$-estimators}
\author{Yiyuan She}
\affil{Department of Statistics\\ Florida State University, Tallahassee, FL 32306}
\date{}
\maketitle

\begin{abstract}
In large-scale modern data analysis,    first-order optimization methods are usually favored to obtain sparse estimators in high dimensions. This paper performs theoretical analysis of  a  class of iterative thresholding based estimators defined in this way. Oracle inequalities are built to show the nearly  minimax rate optimality of such estimators under a new type of  regularity conditions. Moreover, the sequence of iterates is found to be able to approach the statistical truth within the best statistical accuracy  geometrically fast. Our results also reveal different benefits brought by  convex and nonconvex types of shrinkage.
\end{abstract}

\section{Introduction}
\label{sec:intro}
Big data naturally arising in  machine learning, biology,   signal processing, and many other areas, call for the need of scalable optimization in computation. Although for low-dimensional problems, Newton or quasi-Newton methods converge fast and have efficient implementations,  they typically do not scale well    to high dimensional data. In contrast,  \textit{first-order} optimization methods    have recently attracted a great deal of attention from researchers  in statistics,  computer science and engineering.
They      iterate  based on  the gradient (or a subgradient) of the objective function, and have      each iteration step  being cost-effective.
In high dimensional statistics, a first-order  algorithm typically proceeds  in the following manner
\begin{align}
\bsbb^{(t+1)} = \mathcal P \circ (\bsbb^{(t)} - \alpha\nabla l(\bsbb^{(t)})),\label{firstorder}
\end{align}
where $\mathcal P$ is an operator that is easy to compute,   $\nabla l$ denotes the gradient of the loss function $l$, and $\alpha$ gives the stepsize. Such a simple iterative procedure  is   suitable for large-scale optimization, and converges  in arbitrarily high dimensions provided   $\alpha$ is properly small.

  $\mathcal P$ can be motivated from  the perspective of statistical    shrinkage or regularization  and is        necessary  to achieve good   accuracy  when the dimensionality is moderate or high. For example,   a   proximity operator \citep{prox14}  is  associated with a convex penalty function.   
But the problems of interest may not always be convex. Quite often,    $\mathcal P$ is  taken as  a certain  \textbf{thresholding rule} $\Theta$ in  statistical learning, such as   SCAD \citep{Fan01_SCAD}.  The resulting   computation-driven estimators, which we call \textit{$\Theta$-estimators},   are   fixed points of
$
\bsbb= \Theta(\bsbb - \nabla l(\bsbb); \lambda)
$.    To study  the non-asymptotic behavior of $\Theta$-estimators (regardless of  the  sample size and   dimensionality), we   will establish  some  oracle inequalities.

During the last decade, people  have performed       rigorous finite-sample   analysis      of many high-dimensional estimators defined as \textit{globally} optimal solutions to    some convex or    nonconvex problems---see \cite{Bunea}, \cite{ZhangHuang}, \cite{bickel09}, \cite{lounici-2010}, \cite{ZZconcave}, \cite{she2014selectable},  among many others.     $\Theta$-estimators    pose   some new questions. First,  although nicely,  an associated optimization criterion  can be constructed for any given $\Theta$-estimator, the objective may not be convex, and the estimator may not correspond to any functional local (or global) minimum. Second, there are various types of $\Theta$-estimators  due to the  abundant     choices of $\Theta$, but a  comparative study regarding their  statistical performance in high dimensions is lacking in the literature.
Third, $\Theta$-estimators are usually computed in an inexact way on big datasets. Indeed, most practitioners   (have to) terminate        \eqref{firstorder}   before full computational convergence.
These disconnects between    theory and practice when using iterative thresholdings   motivate our work.

The rest of the paper is organized as follows. Section \ref{sec:bg} introduces the $\Theta$-estimators, the associated iterative algorithm--TISP, and some necessary notation. Section \ref{sec:main} presents the main results, including some oracle inequalities, and sequential analysis of the iterates generated by TISP. Section \ref{sec:proofs} provides   proof details.

\section{Background and Notation}
\label{sec:bg}
\subsection{Thresholding functions}
\begin{definition}[Thresholding function]\label{def:threshold}
A thresholding function is a real valued
function  $\Theta(t;\lambda)$ defined for $-\infty<t<\infty$
and $0\le\lambda<\infty$ such that
(i) $\Theta(-t;\lambda)= -\Theta(t;\lambda)$;
(ii) $\Theta(t;\lambda)\le \Theta(t';\lambda)$ for $t\le t'$; 
(iii) $\lim_{t\to\infty} \Theta(t;\lambda)=\infty$;
(iv) $0\le \Theta(t;\lambda)\le t$\ for\ $0\le t<\infty$.
\end{definition}
A vector version of $\Theta$ (still denoted by $\Theta$) is defined componentwise if
either $t$ or $\lambda$ is replaced by a vector. 
From the definition,
\begin{align}\Theta^{-1}(u;\lambda):= \sup\{t:\Theta(t;\lambda)\leq u\}, \forall u > 0
\end{align}   must be monotonically non-decreasing  and so its derivative is defined almost everywhere on $(0, \infty)$. Given $\Theta$, a     critical number   ${\mathcal L}_{\Theta}\le 1$ can be introduced  such that  $\rd\Theta^{-1}(u;\lambda)/\rd u\ge 1- {\mathcal L}_{\Theta}$ for  {almost every} $u\ge 0$, or
\begin{align}
{\mathcal L}_{\Theta} := 1- \essinf\{ \rd \Theta^{-1}(u;\lambda)/\rd u: u \ge 0\},
\end{align}
where    $\essinf$ is the   essential infimum.
For the perhaps most popular soft-thresholding and hard-thresholding functions
$$\Theta_S(t;\lambda)= \mbox{sgn} (t)(|t| - \lambda)_+, \quad \Theta_H(t;\lambda)=t 1_{|t| \ge \lambda},$$
 ${\mathcal L}_{\Theta}$ equals    $0$ and $1$, respectively.

For any  arbitrarily given $\Theta$,  we construct a penalty function $P_{\Theta}(t;\lambda)$   as follows
\begin{align}
P_{\Theta}(t; \lambda)=\int_0^{|t|} (\Theta^{-1}(u;\lambda) - u) \rd u= \int_0^{|t|} (\sup\{s:\Theta(s;\lambda)\leq u\} - u) \rd u  \label{pendef}
\end{align}
for any $t\in \mathbb R$.
This penalty will be used to make a proper objective function for $\Theta$-estimators.

The threshold $\tau(\lambda) :=\Theta^{-1}(0; \lambda)$  may not equal $\lambda$ in general. For ease in notation,  in writing $\Theta(\cdot;\lambda)$, we  always assume that $\lambda$ is the \emph{threshold parameter}, i.e., $\lambda=\tau(\lambda)$, unless otherwise specified. Then an important fact is that given  $\lambda$, any thresholding rule $\Theta$ satisfies $\Theta(t;\lambda) \le \Theta_H(t;\lambda)$, $\forall t\ge 0$, due to property (iv), from which it  follows that
\begin{align} P_{\Theta} (t; \lambda) \ge P_H(t; \lambda),
\end{align}
 where
\begin{align}
P_H(t; \lambda)=\int_0^{|t|} (\Theta_H^{-1}(u;\lambda) - u) \rd u =(-t^2/2+\lambda |t|)1_{|t|<\lambda} +(\lambda^2/2) 1_{|t|\geq \lambda}.\label{PHdef}
\end{align}
In particular,   $P_H(t; \lambda) \le P_0(t; \lambda):= \frac{\lambda^2}{2} 1_{t\neq 0}$ and $P_H(t; \lambda) \le P_1(t; \lambda):= \lambda |t|$.

When $\Theta$ has discontinuities, such as $t= \pm \lambda$ in $\Theta_H(t; \lambda)$, ambiguity  may arise in  definition. To avoid the issue, we  assume the quantity to be thresholded never corresponds to any discontinuity of $\Theta$. This assumption is mild  because practically used thresholding rules  have  few discontinuity points and such discontinuities rarely occur in   real applications.
\subsection{$\Theta$-estimators}\label{introThetaest}
We assume a model
\begin{align}
\bsby = \bsbX \bsbb^* + \bsbeps,
\end{align}
where $\bsbX$ is an ${n\times p}$  design matrix, $\bsby $ is a response vector in $\mathbb R^{n}$, $\bsbb^*$ is the unknown coefficient vector, and $\bsbeps$ is a \textit{sub-Gaussian} random vector with mean zero and scale bounded by $\sigma$,   cf. Definition \ref{def:subgauss} in Section \ref{sec:proofs} for more detail. Then a $\Theta$-estimator $\hat \bsbb$, driven by the computational procedure \eqref{firstorder}, is defined as a solution to the    $\Theta$-equation
\begin{align}
\rho\bsbb = \Theta(\rho \bsbb + \bsbX^T \bsby/\rho - \bsbX^T \bsbX \bsbb/\rho; \lambda), \label{scaledthetaeq}
\end{align}
where  $\rho$, the scaling parameter, does not depend on $\bsbb$. 
Having $\rho$ appropriately large   is crucial to guarantee the convergence of the computational procedure.

 All popularly used penalty functions are associated with thresholdings, such as the $\ell_r$ ($0< r \le 1$), $\ell_2$, SCAD \citep{Fan01_SCAD}, MCP \citep{mplus}, capped $\ell_1$ \citep{Zhang10}, $\ell_0$, elastic net \citep{ZouHas}, Berhu \citep{Owen,s2}, $\ell_0+\ell_2$ \citep{SheTISP}, to name a few.     Table \ref{tab:thetaexs} lists some examples. From  a shrinkage perspective, thresholding rules usually  suffice in statistical learning.

Equation \eqref{scaledthetaeq} can be re-written  in terms of the scaled deign   $\tilde\bsbX  = \bsbX/\rho$ and the corresponding coefficient vector   $\tilde\bsbb = \rho\bsbb $
\begin{align}
\tilde\bsbb = \Theta( \tilde\bsbb + \tilde\bsbX^T \bsby - \tilde\bsbX^T \tilde\bsbX \tilde\bsbb; \lambda). \label{scaledthetaeq2}
\end{align}
 We will show that the     $\lambda$        in the scaled form does not have to adjust for the sample size, which is advantageous   in regularization parameter tuning.

A simple  iterative procedure can be defined based on \eqref{scaledthetaeq} or \eqref{scaledthetaeq2}: \begin{align}\tilde\bsbb^{(t+1)} = \Theta( \tilde\bsbb^{(t)} + \tilde\bsbX^T \bsby - \tilde\bsbX^T \tilde\bsbX \tilde\bsbb^{(t)}; \lambda) , \bsbb^{(t+1)}= \tilde\bsbb^{(t+1)}/\rho, \label{tispdef}\end{align}  which is called   the Thresholding-based Iterative Selection Procedure (\textbf{TISP}) \citep{SheTISP}.
From Theorem 2.1 of \cite{SheGLMTISP}, given an arbitrary $\Theta$, TISP ensures the following function-value descent property when $\rho\ge \frac{\| \bsbX\|_2}{2 - {\mathcal L}_{\Theta}}$:
\begin{align}
f(\bsbb^{(t+1)}; \lambda) \le f(\bsbb^{(t)}; \lambda).
\end{align}
Here,  the energy function (objective function) is constructed as \begin{align}
f(\bsbb;\lambda) &= \frac{1}{2} \| \bsbX \bsbb- \bsby\|_2^2 + \sum_{j=1}^p P(\rho|\beta_j|; \lambda),\label{objfunc}
\end{align}
where the penalty $P$ can be $P_{\Theta}$ as defined  in \eqref{pendef}, or more generally,
\begin{align}
P(t;\lambda)&= P_{\Theta}(t; \lambda)
+ q(t; \lambda),
\end{align}
with $q$ an arbitrary function   satisfying $q(t,\lambda)\ge 0$, $\forall t\in \mathbb R$   and  $q(t;\lambda)=0$ if $t=\Theta(s;\lambda)$ for some $s\in \mathbb R$.
Furthermore, we can show that  when $\rho>{\| \bsbX\|_2}/({2 - {\mathcal L}_{\Theta}})$, any limit point of $\bsbb^{(t)}$ is  necessarily a fixed point of \eqref{scaledthetaeq}, and thus a $\Theta$-estimator. See \cite{SheGLMTISP} for more detail.
Therefore, $f$ is not necessarily unique when $\Theta$  has discontinuities---for example,  penalties like the capped $\ell_1$, $P_0(t;\lambda) = \frac{\lambda^2}{2} 1_{t\ne 0}$ and $P_H$  are all associated with the same  $\Theta_H$.
Because of the  \textit{many-to-one} mapping from penalty functions to thresholding functions, iterating \eqref{firstorder} with a well-designed thresholding rule is perhaps  more convenient than solving a nonconvex penalized optimization problem. Indeed, some penalties (like SCAD) are designed from the thresholding viewpoint.

\begin{table}[!ht]
\centering
\caption{\small{Some examples  of thresholding functions and their associated quantities.} }\label{tab:thetaexs}

\scriptsize{
\addtolength{\tabcolsep}{5pt}
\begin{tabular}{l  l l l}
\hline

\hline

\hline
 & \textbf{soft} & \textbf{ridge} & \textbf{hard}    \\
\hline
$\Theta$ &
$(t - \lambda\mbox{sgn}(t)) 1_{|t|>\lambda}$ & $\frac{t}{1+\eta}$ & $t 1_{|t|>\lambda}$  \\

\hline
$\mathcal L_{\Theta}$ & $0$ & $-\eta$ & $1$   \\
\hline
$P_{\Theta}$ & $\lambda|t|$ & $\frac{\eta}{2} t^2$ & $\begin{cases} -\frac{1}{2} t^2 + \lambda |t|, &\mbox{ if } |t| < {\lambda}\\ \frac{1}{2}{\lambda^2}, &\mbox{ if } |t| \geq {\lambda} \end{cases}$
\\
 $P$ & & & $ \min(\lambda |t|, \frac{\lambda^2}{2})$ (`capped $\ell_1$'), $ \frac{\lambda^2}{2} 1_{t\ne 0} $
 \\
 \hline

 \hline
 & \textbf{elastic net} ($\eta\ge0$) & \textbf{berhu} ($\eta\ge0$) & \textbf{hard-ridge} ($\eta\ge0$)    \\
\hline
$\Theta$ &
$\frac{t-\lambda\mbox{sgn}(t)}{1+\eta}1_{|t|\geq \lambda}$ & $ \begin{cases} 0 & \text{if }
|t|<\lambda\\
t-\lambda \mbox{sgn}(t) & \text{if } \lambda \leq |t|\leq
\lambda+\lambda/\eta\\
\frac{t}{1+\eta} & \text{if } |t|>\lambda+\lambda/\eta
\end{cases}
$ & $\frac{t}{1+\eta} 1_{|t|>\lambda}$  \\

\hline
$\mathcal L_{\Theta}$ & $-\eta$ & $0$ & $1$   \\
\hline
$P_{\Theta}$ & $\lambda|t|+\frac{1}{2}\eta t^2$ & $
\begin{cases}
\lambda|t| & \text{if } |t|\leq \lambda/\eta\\
\frac{\eta t^2}{2}+\frac{\lambda^2}{2\eta} & \text{if } |t|> \lambda/\eta. \end{cases}
$ & $\begin{cases} -\frac{1}{2} t^2 + \lambda |t|, &\mbox{ if } |t| < \frac{\lambda}{1+\eta}\\ \frac{1}{2} \eta t^2 +\frac{1}{2}\frac{\lambda^2}{1+\eta}, &\mbox{ if } |t| \geq \frac{\lambda}{1+\eta}. \end{cases}$
 \\
 $P$ & & & $ \frac{1}{2} \frac{\lambda^2}{1+\eta} 1_{t\ne 0} + \frac{\eta}{2} t^2  $  (`$\ell_0+\ell_2$')\\
\hline

\hline
& \multicolumn{2}{l}{\textbf{scad} ($a> 2$)} & \textbf{mcp} ($\gamma\ge 1$) \\
\hline
$\Theta$  & \multicolumn{2}{l}{ $\begin{cases} 0, &\mbox{if } |t| \leq \lambda\\ t-\lambda\,\mbox{sgn}(t), &\mbox{if } \lambda< |t| \leq 2\lambda\\ \frac{(a-1)t-a\lambda\,\mbox{sgn}(t)}{a-2}, & \mbox{if } 2\lambda<|t| \leq a\lambda \\ t, & \mbox{if } |t| > a\lambda \end{cases}
$ }  & $\begin{cases}0, & \mbox{ if } |t|< \lambda\\ \frac{t-\lambda\mbox{sgn}(t)}{1-1/\gamma}, & \mbox{ if } \lambda \leq   |t| < \gamma \lambda\\
t, & \mbox{ if }     |t|  \geq \gamma \lambda
\end{cases}$
\\

\hline
$\mathcal L_{\Theta}$ & \multicolumn{2}{l}{ $1/(a-1)$ } &  $1/\gamma$  \\
\hline
$P_{\Theta}$  & \multicolumn{2}{l}{ $\frac{{\mathrm d} P}{{\mathrm d} t} = \begin{cases} \lambda\,\mbox{sgn}(t), &\mbox{if } |t|\leq\lambda\\ \frac{a\lambda\,\mbox{sgn}(t)-t}{a-1}, &\mbox{if } \lambda<|t|\leq a\lambda\\ 0, & \mbox{if } |t| > a \lambda\end{cases}
$ } & $\begin{cases} -\frac{t^2}{2\gamma}  + \lambda |t|,&\mbox{if } |t| < \gamma{\lambda}\\ \frac{\gamma\lambda^2}{2},&\mbox{if } |t| \geq \gamma{\lambda} \end{cases}=\frac{1}{\gamma}P_H(t; \gamma \lambda) $
 \\
\hline

\hline
& \multicolumn{3}{l}{$\bf l_r$ ($0<r<1$, $\zeta\ge 0$)}   \\
\hline
$\Theta$  & \multicolumn{3}{l}{ $\begin{cases}0, \mbox{ if } |t| \leq \zeta^{1/(2-r)} (2-r)(2-2r)^{(r - 1)/(2-r)}  \\
    \mbox{sgn}(t)\max\{\zeta^{1/(2-r)} [r(1-r)]^{1/(2-r)} \le \theta\le |t|: \theta+\zeta r \theta^{r-1} = |t|\},  \mbox{ otherwise}. (\mbox{The set is a singleton.}) \end{cases}
$ }  \\

\hline
$\mathcal L_{\Theta}$ & \multicolumn{3}{l}{ $1$  }   \\
\hline
$P$  & \multicolumn{3}{l}{ $\zeta |t|^r$  }
 \\

\hline

\hline

\hline
\end{tabular}
}
\end{table}

The following theorem shows that     the set of  $\Theta$-estimators include all  locally optimal solutions of $\frac{1}{2} \| \bsbX \bsbb- \bsby\|_2^2 + \sum_{j=1}^p P_{\Theta}(|\beta_j|; \lambda)=:f_{\Theta}(\bsbb)$.

\begin{theorem}\label{localopt}
 Let $\hat\bsbb$ be a  local minimum point (or a coordinate-wise minimum point) of $f_{\Theta}(\cdot)$. If $\Theta$ is continuous at $\hat \bsbb + \bsbX^T \bsby - \bsbX^T \bsbX\hat \bsbb$,
 $\hat\bsbb$ must satisfy $\bsbb = \Theta( \bsbb + \bsbX^T \bsby - \bsbX^T \bsbX \bsbb; \lambda)$. 
\end{theorem}
The converse is not necessarily true. Namely,  $\Theta$-estimators may not guarantee functional     local optimality, let alone  global optimality.
This raises   difficulties in statistical analysis. We will give a novel and unified treatment which can yield nearly optimal error rate for various thresholdings.

\section{Main Results}
\label{sec:main}
To address the problems   in arbitrary  dimensions (with possibly large $p$ and/or $n$), we aim to establish non-asymptotic {oracle inequalities} \citep{donoho1994}.    For any $\bsbb= [ \beta_1, \ldots, \beta_p]^T$, define
\begin{align}\mathcal J(\bsbb)= \{j: \beta_j \neq {0}\}, \qquad  J(\bsbb)=|\mathcal J(\bsbb)|=\|\bsbb\|_{0}.
\end{align}
Recall $P_1(t; \lambda) = \lambda |t|$,  $P_0(t; \lambda) = \frac{\lambda^2}{2} 1_{t\neq 0}$, $P_H(t; \lambda)= (-t^2/2+\lambda |t|)1_{|t|<\lambda} +(\lambda^2/2) 1_{|t|\geq \lambda}.
$ For convenience, we use   $P_{1} (\bsbb; \lambda)$ to denote $\lambda \| \bsbb\|_{1}$ when there is no ambiguity.   $P_{0}(\bsbb; \lambda)$ and $P_{H}(\bsbb; \lambda)$  are used    similarly.
We denote by     $\lesssim$    an inequality that holds up to a multiplicative  constant.

Unless otherwise specified, we  study     \emph{scaled}  $\Theta$-estimators satisfying  equation   \eqref{scaledthetaeq2}, where $\tilde\bsbb = \rho\bsbb $,  $\tilde \bsbX  = \bsbX/\rho$, and   $\rho\ge {\| \bsbX\|_2}$ (and so $\|\tilde \bsbX\|_2\le 1$). By abuse of notation,  we still write  $\bsbb$ for $\tilde \bsbb$, and   $\bsbX$ for $\tilde\bsbX$.
As mentioned previously, we always   assume that $\Theta$ is continuous at $\hat\bsbb + \bsbX^T \bsby - \bsbX^T \bsbX \hat\bsbb$ in Sections \ref{subsec:PTtype} \&\  \ref{subsec:P0type}; similarly,     Section \ref{subsec:seq}  assumes that $\Theta$ is continuous at $\bsbb^{(t)} + \bsbX^T \bsby - \bsbX^T \bsbX \bsbb^{(t)}$.

The past works on the lasso show that  a certain incoherence requirement  must be assumed   to obtain sharp error rates. In most theorems, we also need to make similar assumptions to prevent the design matrix  from being too collinear. We will state  a new type of regularity conditions, which are called \emph{comparison regularity conditions},   under which   oracle inequalities and sequential statistical error bounds can be  obtained for any $\Theta$.
\subsection{$P_{\Theta}$-type oracle inequalities under ${\mathcal  R}_0$}
\label{subsec:PTtype}
In this subsection, we  use   $P_{\Theta}$ to make a bound of the prediction error of  $\Theta$-estimators. Our \ regularity condition is stated  as follows.\\

\noindent{\textsc{Assumption ${\mathcal  R}_0(\delta, \vartheta, K, \bsbb, \lambda)$}} Given $\bsbX$, $\Theta$,  $\bsbb$, $\lambda$, there exist  $\delta > 0$, $\vartheta>0$, $K\ge0$ such that the following inequality holds for any $\bsbb'\in \mathbb R^p$
\begin{align}
\begin{split}
    &{\vartheta}  P_{H}(\bsbb' - \bsbb; {\lambda}) + \frac{\mathcal L_{\Theta}}{2} \|\bsbb' - \bsbb\|_2^2
\\\le \ &  \frac{2-\delta}{2} \| \bsbX  (\bsbb' - \bsbb)\|_2^2 + P_{\Theta}( \bsbb' ; \lambda) + K  P_{\Theta}(\bsbb; \lambda).
\end{split}\label{reg0}
\end{align}

Roughly, \eqref{reg0} means that       $2 \| \bsbX  (\bsbb' - \bsbb)\|_2^2$ can   dominate ${\mathcal L_{\Theta}} \|\bsbb' - \bsbb\|_2^2$  with the help from $P_{\Theta}( \bsbb' ; \lambda)$ and $K  P_{\Theta}(\bsbb; \lambda)$ for some $K>0$. 

\begin{theorem}\label{th_oracle}
Let   $\hat\bsbb$ be  any $\Theta$-estimator satisfying $ \bsbb = \Theta(  \bsbb +  \bsbX^T \bsby -  \bsbX^T  \bsbX \bsbb; \lambda)$ with  $\lambda =A\sigma \sqrt{\log (ep)}$ and $A$ a constant. Then  for any sufficiently large $A$, the following  oracle  inequality holds for    $\bsbb\in \mathbb R^{p}$
\begin{align}
 \EE [ \| \bsbX \hat \bsbb - \bsbX \bsbb^* \|_2^2 ]\lesssim    \|\bsbX  \bsbb - \bsbX \bsbb^* \|_2^2 +
P_{\Theta}(\bsbb; \lambda) +\sigma^{2},\label{genoracle}
\end{align}
provided   ${\mathcal  R}_0(\delta, \vartheta, K, \bsbb, \lambda)$ is satisfied for some constants $\delta>0, \vartheta>0$, $K\ge0$.
\end{theorem}

 Theorem \ref{th_oracle}    is applicable to any $\Theta$. Let's   examine two specific cases.
 First, consider    $\mathcal L_{\Theta}\le 0$, which indicates that $P_{\Theta}$ is convex. Because $P_H \le P_{\Theta}$ and $P_H$ is sub-additive: $P_{H}(t+s) \le P_{H}(t) + P_{H}(s)$  due to its concavity \citep{ZZconcave}, ${\mathcal  R}_0(\delta, \vartheta, K, \bsbb, \lambda)$ is \textit{always} satisfied (for any  $\delta \le 2$, $0<\vartheta\le 1$,    $K\ge \vartheta$).

\begin{corollary} \label{conv-regfree-oracle}
Suppose $\Theta$ satisfies  $\mathcal L_{\Theta}\le 0$. Then, \eqref{genoracle} holds for all corresponding $\Theta$-estimators, without requiring any regularity condition. \end{corollary}

In the case of hard-thresholding or SCAD thresholding, $P_{\Theta}(\bsbb;\lambda)$ does not depend on the magnitude of $\bsbb$, and we can get   a finite complexity rate in the oracle inequality. Also, $\mathcal R_0$ can be slightly relaxed,  by replacing   $K  P_{\Theta}(\bsbb; \lambda)$ with $K  P_0(\bsbb; \lambda)$ in \eqref{reg0}. We denote   the modified version     by ${\mathcal  R}_0'(\delta, \vartheta, K, \bsbb, \lambda)$.
\begin{corollary}\label{l0like-oracle}
  Suppose that $\Theta$ corresponds to  a bounded nonconvex penalty satisfying   $P_{\Theta}(t;\lambda) \leq C \lambda^2$, $\forall t \in \mathbb R$, for some constant $C>0$. 
 Then in the setting of Theorem \ref{th_oracle},   \begin{align}
 \EE [ \| \bsbX \hat \bsbb - \bsbX \bsbb^* \|_2^2 ]\lesssim    \|\bsbX  \bsbb - \bsbX \bsbb^* \|_2^2 +
\sigma^2 J(\bsbb) \log ( e p) +\sigma^{2}, \label{rateoracle}
\end{align}
provided  ${\mathcal  R}_0'(\delta, \vartheta, K, \bsbb, \lambda)$  is satisfied for some constants $\delta>0, \vartheta>0$, $K\ge0$.
\end{corollary}

\begin{remark} \upshape \label{rem2} The right-hand side of the oracle inequalities involves a \textit{bias} term $ \| \bsbX \bsbb - \bsbX \bsbb^*\|_2^2$ and a \textit{complexity}   term $P_{\Theta}(\bsbb; \lambda)$.
 Letting     $\bsbb = \bsbb^*$ in, say, \eqref{genoracle},   the bias    vanishes,  and
we obtain  a prediction error bound   of the order  $\sigma^2 J^* \log (ep) $  (omitting  constant factors), where $J^*$ denotes  the number of nonzero components in  $\bsbb^*$.
On the other hand,  the existence of the bias term ensures the  applicability of  our results to  approximately  sparse  signals. For example,  when $\bsbb^*$ has many small but nonzero components,
we can use a reference $\bsbb$ with a  much smaller   support than $\mathcal J(\bsbb^*)$ to get a lower error bound, as a benefit from the bias-variance tradeoff.
\end{remark}

\begin{remark} \upshape \label{rem0}
When $\mathcal R_0$ holds with $\delta>1$,  the proof of Theorem \ref{th_oracle}  shows that the multiplicative constant for   $\|\bsbX  \bsbb - \bsbX \bsbb^* \|_2^2$  can be as small as $1$.    The corresponding oracle inequalities are    called `sharp' in some works \citep{Kol11}.  This also applies to Theorem \ref{th_oracle-l1type}. 
    Our proof scheme can also deliver high-probability form  results, without   requiring an upper bound of $\|\bsbX\|_2$. 
\end{remark}



\begin{remark} \upshape \label{rem3} Corollary \ref{l0like-oracle} applies to all  ``hard-thresholding like'' $\Theta$, because when $\Theta(t;\lambda) =t$ for $|t|>c \lambda$,   $P_{\Theta}(t;\lambda) \le c^2 \lambda^2$. It is worth mentioning that   the  error rate of $\sigma^2 J^* \log (ep)$  cannot be significantly improved in a minimax sense.   In fact, under the Gaussian noise contamination and some regularity conditions,   there exist   constants $C$, $c>0$  such that
$
    \inf_{\check \bsbb}\,\sup_{\bsbb^*: J(\bsbb^*)\leq J   }\mathbb{E}[\| \bsbX (\check{\bsbb}-\bsbb^*)\|_2^2)/(C P_o(J))] \geq c>0, 
$
where   $\check\bsbb$ denotes an arbitrary estimator of $\bsbb^*$  and
$
P_o(J)= \sigma^2\{ J + J\log(ep/J)\}.
$ See, e.g., 
 \cite{lounici-2010} for a proof.
     The bound in \eqref{rateoracle}   achieves  the minimax  optimal rate up to a mild logarithm factor  for any $n$ and $p$.
\end{remark} %



\subsection{$P_0$-type oracle inequalities under ${\mathcal  R}_1$}
\label{subsec:P0type}
This part uses $P_0$ instead of $P_{\Theta}$ to make an oracle bound. We will show that under another type of comparison regularity conditions,  \textit{all} thresholdings  can  attain the essentially optimal error rate given in Corollary \ref{l0like-oracle}. We will also show that in the   case of soft-thresholding, our condition is  more relaxed  than many other assumptions   in the literature.  \\


\noindent{\textsc{Assumption ${\mathcal  R}_1(\delta, \vartheta, K, \bsbb, \lambda)$}} Given $\bsbX$, $\Theta$,  $\bsbb$, $\lambda$, there exist  $\delta > 0$, $\vartheta>0$, $K\ge0$ such that the following inequality holds for any $\bsbb'\in \mathbb R^p$
\begin{align}
\begin{split}
&   \vartheta  P_{H}(\bsbb' - \bsbb; {\lambda}) + \frac{\mathcal L_{\Theta}}{2} \|\bsbb' - \bsbb\|_2^2 + P_{\Theta}(\bsbb; \lambda)  \\
\le \ & \frac{2-\delta}{2} \| \bsbX  (\bsbb' - \bsbb)\|_2^2  +  P_{\Theta}( \bsbb' ; \lambda)+ K \lambda^2 {J}(\bsbb).
\end{split}\label{reg1}
\end{align}

\begin{theorem}\label{th_oracle-l1type}
Let     $\hat\bsbb$ be  a $\Theta$-estimator and $\lambda =A\sigma \sqrt{\log (ep)}$  with $A$ a sufficiently large constant. Then   $\EE [ \| \bsbX \hat \bsbb - \bsbX \bsbb^* \|_2^2 ]\lesssim    \|\bsbX  \bsbb - \bsbX \bsbb^* \|_2^2 +
\lambda^2 J(\bsbb)  +\sigma^{2}$ holds for  any $\bsbb\in \mathbb R^{p}$
if ${\mathcal  R}_1(\delta, \vartheta, K, \bsbb, \lambda)$ is satisfied for some constants $\delta>0, \vartheta>0$, $K\ge0$.
 \end{theorem}

\begin{remark} \upshape \label{rem6}
 Some fusion thresholdings, like those associated with elastic net, Berhu  and Hard-Ridge (cf. Table \ref{tab:thetaexs}), involve an additional $\ell_2$ shrinkage. In the \ situation,   the complexity term in the oracle inequality should involve both $J(\bsbb)$ and $\| \bsbb\|_2^2$.
We can modify our regularity   conditions to obtain such $\ell_0+\ell_2$  bounds using the same proof scheme. The details are  however not reported in this paper.
In addition, our results can be   extended to   $\Theta$-estimators with a stepsize parameter. Given $\lambda>0$ and $0<\alpha \le 1$, suppose  $\lambda_{\alpha}$ is introduced   such that $\alpha P_{\Theta}(t;\lambda) = P_{\Theta}(t; \lambda_{\alpha})$  for any $\ t$. Then, for any $\hat \bsbb$ as a fixed point of
$
\bsbb = \Theta(\bsbb - \alpha \bsbX^T \bsbX \bsbb + \alpha \bsbX^T\bsby; \lambda_{\alpha}), 
$
an analogous result can be obtained (the only change is that   ${\mathcal L}_{\Theta}$ is replaced by ${\mathcal L}_{\Theta}/\alpha$).
\end{remark} 

To give some more intuitive regularity conditions, we suppose $P_{\Theta}$ is concave on $[0, \infty)$. Examples include     $\ell_r$ ($0\le r \le 1$),  MCP,   SCAD, and so on. The concavity    implies    $P_{\Theta}(t+s) \le P_{\Theta}(t) + P_{\Theta}(s)$, and so      $P_{\Theta}(\bsbb_{\mathcal J}';\lambda) - P_{\Theta}(\bsbb_{\mathcal J};\lambda)\le P_{\Theta}((\bsbb' - \bsbb)_{\mathcal J};\lambda)$ and  $P_{\Theta}(\bsbb_{\mathcal J^c}';\lambda) = P_{\Theta}((\bsbb' - \bsbb)_{\mathcal J^c};\lambda)$, where $\mathcal J^c$ is the complement of $\mathcal J$ and   $\bsbb_{\mathcal J}$ is the subvector of $\bsbb$ indexed by $\mathcal J$. Then $\mathcal R_1$ is implied by
${\mathcal  R}_1'$   below for  given $\mathcal J = \mathcal J(\bsbb)$. \\

\noindent{\textsc{Assumption ${\mathcal  R}_1'(\delta, \vartheta, K, \mathcal J, \lambda)$}} Given $\bsbX$, $\Theta$,  $\mathcal J$, $\lambda$, there exist  $\delta > 0$, $\vartheta>0$, $K\ge0$ such that  for any $\bsbDelta\in \mathbb R^p$,
 \begin{align}
 \begin{split}
&    {  P_{\Theta}(\bsbDelta_{\mathcal J}; \lambda)+\vartheta}  P_{H}(\bsbDelta_{\mathcal J}; {\lambda}) + \frac{\mathcal L_{\Theta}}{2} \|\bsbDelta\|_2^2   \\
\le \ & \frac{2-\delta}{2} \| \bsbX  \bsbDelta\|_2^2  + K \lambda^2 {J}  +  P_{\Theta}( \bsbDelta_{\mathcal J^c} ; \lambda) - {\vartheta}  P_{H}(\bsbDelta_{\mathcal J^c}; {\lambda}),
\end{split}\label{reg1-2}
\end{align}
or
\begin{align}
\begin{split}
   (1+ \vartheta)   P_{\Theta}(\bsbDelta_{\mathcal J}; \lambda)+ \frac{\mathcal L_{\Theta}}{2} \|\bsbDelta\|_2^2
\le   \frac{2-\delta}{2} \| \bsbX  \bsbDelta\|_2^2  + K \lambda^2 {J}  +  (1-\vartheta) P_{\Theta}( \bsbDelta_{\mathcal J^c} ; \lambda).
\end{split} \label{reg1-3}
\end{align}

When $\Theta$ is the soft-thresholding,  it is easy to verify that    a sufficient condition  for \eqref{reg1-3}  is
\begin{align}
   (1+ \vartheta)   \|\bsbDelta_{\mathcal J}\|_1
\le  K \sqrt{J} \|\bsbX\bsbDelta\|_2 +  \| \bsbDelta_{\mathcal J^c} \|_1, \label{reg1-l1}
\end{align}
for some   $\vartheta> 0$ and $K\ge 0$. \eqref{reg1-l1} has a  simper form than $\mathcal R_1$. In the following, we give the definitions of    the RE  and the compatibility condition \citep{bickel09,van2009conditions} to make a comparison to \eqref{reg1-l1}.  \\

\noindent{\textsc{Assumption ${\mathcal  RE}(\kappa_{RE}, \vartheta_{RE}, \mathcal J)$}}. \
Given $\mathcal J\subset [p]$, we say that  $\bsbX\in {\mathbb R}^{n\times p}$ satisfies $\mathcal RE(\kappa_{RE}, \vartheta_{RE}, \mathcal J)$, if for positive numbers $\kappa_{RE}$, $\vartheta_{RE}>0$,
 \begin{align}
{J} \| \bsbX \bsbDelta \|_2^2   \ge\kappa_{RE}\| \bsbDelta_{\mathcal J}\|_1^2, \label{conditionComp}
 \end{align}
 or more restrictively,
 \begin{align}
\| \bsbX \bsbDelta \|_2^2   \ge\kappa_{RE}  \| \bsbDelta_{\mathcal J}\|_2^2,\label{conditionRE}
 \end{align}
  for all $\bsbDelta\in \mathbb R^{p}$ satisfying
\begin{align}
 (1+\vartheta_{RE})\| \bsbDelta_{\mathcal J}\|_1  \geq \| \bsbDelta_{\mathcal J^c}\|_1. \label{conecond}
\end{align}
\\

Assume  ${\mathcal  RE}(\kappa_{RE}, \vartheta_{RE}, \mathcal J)$
holds. When $(1+ \vartheta_{RE})   \|\bsbDelta_{\mathcal J}\|_1 \le  \| \bsbDelta_{\mathcal J^c} \|_1$, \eqref{reg1-l1} holds trivially with $\vartheta=\vartheta_{RE}$; otherwise,  \eqref{conditionComp} indicates  $(1+ \vartheta)   \|\bsbDelta_{\mathcal J}\|_1 \le  K\sqrt{J} \| \bsbX \bsbDelta \|_2 $ with  $K = (1+\vartheta_{RE})/\sqrt{\kappa_{RE}}$.
So intuitively, we have the following relationship:
$$
\eqref{conditionRE} +\eqref{conecond} \Rightarrow \eqref{conditionComp}+\eqref{conecond} \Rightarrow \eqref{reg1-l1} \Rightarrow  \eqref{reg1-3} \Rightarrow  \eqref{reg1-2} \Rightarrow  \eqref{reg1}.
$$
  In particular,    ${\mathcal  R}_1$ is  less demanding than RE. \\

Next, let's  compare the regularity conditions required by $\Theta_S$ and $\Theta_H$ to achieve the nearly optimal  error rate. Recall    ${\mathcal  R}_1(\delta, \vartheta, K, \bsbb, \lambda)$   and
 ${\mathcal  R}_0'(\delta, \vartheta, K, \bsbb, \lambda)$ in   Theorem \ref{th_oracle-l1type} and Corollary \ref{l0like-oracle}, respectively
\begin{align*}
   \vartheta  P_{H}(\bsbb' - \bsbb; {\lambda}) + \lambda \|\bsbb\|_1
&\le  \frac{2-\delta}{2} \| \bsbX  (\bsbb' - \bsbb)\|_2^2  + \lambda \|\bsbb'\|_1 + K \lambda^2 {J},
\\    {\vartheta}  P_{H}(\bsbb' - \bsbb; {\lambda}) + \frac{1}{2} \|\bsbb' - \bsbb\|_2^2 &\le  \frac{2-\delta}{2} \| \bsbX  (\bsbb' - \bsbb)\|_2^2 + P_{H}( \bsbb' ; \lambda) + K \lambda^2 J. 
\end{align*}
  ${\mathcal  R}_0'(\delta, \vartheta, K, \bsbb, \lambda)$ implies  ${\mathcal  R}_1(\delta, \vartheta, K+1, \bsbb, \lambda)$. Indeed,  for
$\bsbDelta = \bsbb' - \bsbb$,
\begin{align*}
\lambda \|\bsbb\|_1  - \lambda \|\bsbb'\|_1  &\le  \lambda \| \bsbDelta_{\mathcal J}\|_1  - \lambda \| \bsbb_{\mathcal J^c}'\|_1 \\
&\le \frac{1}{2} \lambda^2 J + \frac{1}{2} \| \bsbDelta_{\mathcal J}\|_2^2 - P_{H}(\bsbb_{\mathcal J^c}'; \lambda)\\
& \le \frac{1}{2} \lambda^2 J + \frac{1}{2} \| \bsbDelta_{\mathcal J}\|_2^2 - P_{H}(\bsbb'; \lambda) + P_{H}(\bsbb_{\mathcal J}'; \lambda)\\
& \le \frac{1}{2} \lambda^2 J + \frac{1}{2} \| \bsbDelta_{\mathcal J}\|_2^2 - P_{H}(\bsbb'; \lambda) + P_{0}(\bsbb_{\mathcal J}'; \lambda) \\
& \le   \lambda^2 J + \frac{1}{2} \| \bsbDelta_{\mathcal J}\|_2^2 - P_{H}(\bsbb'; \lambda).
\end{align*}

On the other hand, Corollary \ref{l0like-oracle} studies when     \emph{all} $\Theta_H$-estimators have the optimal performance guarantee, while  practically, one may initialize      \eqref{tispdef} with a carefully chosen  starting point.
\begin{theorem} \label{globaloracle}
Given any $\Theta$, there exists a $\Theta$-estimator (which minimizes \eqref{objfunc}) such that \eqref{genoracle}   holds
without requiring any  regularity condition. In particular, if  $\Theta$ corresponds to  a bounded nonconvex penalty as described in Corollary \ref{l0like-oracle}, then there exists a $\Theta$-estimator such that \eqref{rateoracle} holds free of regularity conditions.
\end{theorem}

  Theorem \ref{globaloracle} does {\textit{not}} place any requirement on $\bsbX$. So it seems that     applying $\Theta_H$  may have some further  advantages in practice. (How to {efficiently} pick a  $\Theta_H$-estimator  to completely remove all regularity conditions  is however beyond the  the scope of the current  paper. For a possible idea of relaxing the conditions, see       Remark \ref{rem8}.)
\\


Finally, we  make a discussion of the scaling parameter  $\rho$. Our results so far  are obtained   after performing    $\bsbX \leftarrow \bsbX /\rho$  with $\rho\ge \| \bsbX\|_2$. The prediction error  is invariant to the transformation. But it   affects the regularity conditions.

 Seen from  \eqref{scaledthetaeq}, $1/\rho^2$ is related to  the stepsize   $\alpha$ appearing in \eqref{firstorder}, also known as the \textit{learning rate} in the machine learning literature.  From  the computational results in Section \ref{introThetaest}, $\rho$  must be     large enough to guarantee   TISP is convergent. The larger the value of $\rho$ is, the smaller the stepsize is (and so the slower the convergence is). Based on   the  machine learning literature,  slow  learning rates are always recommended   when training  a nonconvex   learner   (e.g.,     artificial neural networks).
Perhaps interestingly, in addition to  computational efficiency reasons, all our    statistical analyses caution  against using an extremely  large scaling  when $\mathcal L_{\Theta} > 0$.   For example,    ${\mathcal  R}_0'(\delta, \vartheta, K, \bsbb, \lambda)$ for an unscaled   $\bsbX$ reads  ${\vartheta}  P_{H}(\rho(\bsbb' - \bsbb); {\lambda}) + {\rho^2}  \|\bsbb' - \bsbb\|_2^2/2 \le  ({2-\delta})  \| \bsbX  (\bsbb' - \bsbb)\|_2^2/2 + P_{H}( \rho\bsbb' ; \lambda) + K \lambda^2 J$, which     becomes difficult to hold when $\rho$ is very large. This makes    the statistical error bound   break down easily. Therefore,  a good idea is to have       $\rho$ just appropriately  large (mildly greater than $\| \bsbX\|_2$). The sequential analysis   of the iterates in the next part also supports the point.

\subsection{Sequential Algorithmic Analysis}\label{subsec:seq}
We perform   statistical error analysis of the sequence of iterates defined by TISP:
$
\bsbb^{(t+1)} = \Theta( \bsbb^{(t)} + \bsbX^T \bsby - \bsbX^T \bsbX \bsbb^{(t)}; \lambda),
$
where  $\|\bsbX\|_2\le 1$ and   $\bsbb^{(0)}$
is the starting point. The study   is motivated from the fact that   in large-scale applications, $\Theta$-estimators are seldom computed exactly. Indeed, why bother  to run   TISP till computational convergence?  How does  the statistical accuracy    improve  (or deteriorate)  at $t$ increases?   Lately, there are  some   key   advances on the    topic. For example,
 \cite{agarwal2012} showed  that for {convex} problems (not necessarily strongly  convex),       proximal gradient algorithms can be    geometrically fast to approach a globally optimal solution $\hat\bsbb$ within the desired statistical precision,  under a set of    conditions. We however care about the     statistical error between $\bsbb^{(t)}$ and the genuine    $\bsbb^*$ in this work.

    We will introduce two comparison   regularity conditions (analogous to $\mathcal R_0$ and $\mathcal R_1$) to present  both   $P_{\Theta}$-type and  $P_0$-type error bounds. Hereinafter,    denote   $(\bsbb^T \bsbA \bsbb)^{1/2}$ by $\|\bsbb\|_{\bsbA}$, where $\bsbA$ is a positive semi-definite matrix.
\\

\noindent{\textsc{Assumption ${\mathcal  S}_0(\delta, \vartheta, K, \bsbb, \bsbb', \lambda)$}} Given $\bsbX$, $\Theta$,  $\bsbb$, $\bsbb'$, $\lambda$, there exist  $\delta>0$, $\vartheta>0$, $K\ge0$ such that the following inequality holds 
\begin{align}
\begin{split}
&    {\vartheta}  P_{H}(\bsbb' - \bsbb; {\lambda}) + \frac{\mathcal L_{\Theta}+\delta}{2} \|\bsbb' - \bsbb\|_2^2 \\
\le \ & \| \bsbX  (\bsbb' - \bsbb)\|_2^2 + P_{\Theta}( \bsbb' ; \lambda) + K  P_{\Theta}(\bsbb; \lambda).
\end{split}\label{regS0}
\end{align}

\noindent{\textsc{Assumption ${\mathcal  S}_1(\delta, \vartheta, K, \bsbb, \bsbb', \lambda)$}} Given $\bsbX$, $\Theta$,  $\bsbb$, $\bsbb'$, $\lambda$, there exist  $\delta > 0$, $\vartheta>0$, $K\ge0$ such that the following inequality holds  
\begin{align}
\begin{split}
&   \vartheta  P_{H}(\bsbb' - \bsbb; {\lambda}) + \frac{\mathcal L_{\Theta}+\delta}{2} \|\bsbb' - \bsbb\|_2^2 + P_{\Theta}(\bsbb; \lambda)  \\
\le \ & \| \bsbX  (\bsbb' - \bsbb)\|_2^2  +  P_{\Theta}( \bsbb' ; \lambda)+ (K+1) \lambda^2 {J}(\bsbb).
\end{split}\label{regS1}
\end{align}
\\

  \eqref{regS0} and \eqref{regS1} require  a bit more than \eqref{reg0} and \eqref{reg1}, respectively, due to  $\| \bsbX\|_2\le 1$.  The  theorem and the corollary below perform sequential analysis of the iterates and     reveal  the explicit roles of $\delta, \vartheta, K$  (which   can   often be treated as constants).

\begin{theorem}\label{th_seq}
Suppose     ${\mathcal  S}_0(\delta, \vartheta, K, \bsbb^*,\bsbb^{(t+1)}, \lambda)$   is satisfied for some   $\delta>0,  \vartheta>0$, $K\ge0$, then for  $\lambda =A\sigma \sqrt{\log (ep)}/\sqrt{(\delta \wedge \vartheta)\vartheta}$ with $A$ sufficiently large,
 the following error bound  holds with probability at least $1 - C p ^{-c A^2}$:
\begin{align}
 \frac{1 + \delta}{2}\|  \bsbb^{(t+1)} - \bsbb^* \|_{(\bsbI - \bsbX^T \bsbX)}^2 &\le   \frac{1}{2} \|  \bsbb^{(t)} - \bsbb^* \|_{(\bsbI - \bsbX^T \bsbX)}^2 +(K+1)
P_{\Theta}(\bsbb^*; \lambda),\label{seqBnd0}
\end{align}
where $C, c$ are universal positive constants.

Similarly, under the same choice of regularity parameter,     if  ${\mathcal  S}_1(\delta, \vartheta, K, \bsbb^*, \bsbb^{(t)}, \lambda)$ is satisfied for some   $\delta >0, \vartheta>0$, $K\ge0$,    \eqref{seqBnd1} is true with probability at least $1 - C p ^{-cA^2}$:
 \begin{align}
 \frac{1 + \delta}{2}\|  \bsbb^{(t+1)} - \bsbb^* \|_{(\bsbI - \bsbX^T \bsbX)}^2 &\le   \frac{1}{2} \|  \bsbb^{(t)} - \bsbb^* \|_{(\bsbI - \bsbX^T \bsbX)}^2 +
 (K +1)\lambda^2 J^* .\label{seqBnd1}
\end{align}
\end{theorem}

\begin{corollary} \label{corrseq}
In the setting of Theorem \ref{th_seq},  for  any initial point $\bsbb^{(0)}\in \mathbb R^{p}$,  we have
\begin{align} \|  \bsbb^{(t)} - \bsbb^* \|_{(\bsbI - \bsbX^T \bsbX)}^2 &  \le   \kappa^t\|  \bsbb^{(0)} - \bsbb^* \|_{(\bsbI - \bsbX^T \bsbX)}^2  + \frac{\kappa }{1 - \kappa} K'  P_{\Theta}(\bsbb^*; \lambda), \label{recurL0}\\
\|  \bsbb^{(t)} - \bsbb^* \|_{(\bsbI - \bsbX^T \bsbX)}^2  &\le   \kappa^t\|  \bsbb^{(0)} - \bsbb^* \|_{(\bsbI - \bsbX^T \bsbX)}^2  + \frac{\kappa }{1 - \kappa} K'  \lambda^2 J^*, \label{recurL1}\end{align}
  under  ${\mathcal  S}_0(\delta, \vartheta, K, \bsbb^*,\bsbb^{(s)}, \lambda)$ and  ${\mathcal  S}_1(\delta, \vartheta, K, \bsbb^*, \bsbb^{(s)}, \lambda),   0\le s \le t-1$, respectively, with  probability at least $1 - C p ^{-cA^2}$. Here,   $\kappa = 1/(1+\delta)$,   $K' = 2  (K+1)$.
\end{corollary}

\begin{remark} \upshape \label{rem7}   We can get some   sufficient conditions for $\mathcal S_0$ and $\mathcal S_1$, similar to the  discussions made in Section \ref{subsec:P0type}.  When $\|\bsbX\|_2$ is strictly less than $1$,  \eqref{regS0} can be relaxed to  ${\vartheta}  P_{H}(\bsbb' - \bsbb; {\lambda}) + ({\mathcal L_{\Theta}+\delta}) \|\bsbb' - \bsbb\|_2^2/{2}
\le  ({2+\delta }) \| \bsbX  (\bsbb' - \bsbb)\|_2^2/{2} + P_{\Theta}( \bsbb' ; \lambda) + K  P_{\Theta}(\bsbb; \lambda)$ for some $\delta>0$.
The proof in Section \ref{app:subsec:proofth5} also gives expectation-form results,    with an   additional additive term    $C \sigma^2 / (\delta \wedge \vartheta)$  in the upper bounds.
Similar to  Remark \ref{rem6},    we can also study    $\Theta$-iterates with stepsize   $\alpha$, in which case        the weighting matrix in \eqref{seqBnd0}-\eqref{recurL1}      changes   from $\bsbI-\bsbX^T\bsbX$ to $\bsbI/\alpha-\bsbX^T\bsbX$, and the factor $ ({\mathcal L_{\Theta}+\delta})/{2} $ in \eqref{regS0} and \eqref{regS1} is replaced by     $ {(\mathcal L_{\Theta}+\delta)}/{(2\alpha)} $. 
\end{remark}

\begin{remark} \upshape \label{rem8} Theorem \ref{th_seq} still applies when  $\delta, \vartheta, K$ and $\lambda$ are dependent on $t$. For example, if we use a varying threshold sequence, i.e.,  $
\bsbb^{(t+1)} = \Theta( \bsbb^{(t)} + \bsbX^T \bsby - \bsbX^T \bsbX \bsbb^{(t)}; \lambda^{(t)}),
$
then \eqref{recurL1} becomes
$$\|  \bsbb^{(t)} - \bsbb^* \|_{(\bsbI - \bsbX^T \bsbX)}^2  \le   \kappa^t\|  \bsbb^{(0)} - \bsbb^* \|_{(\bsbI - \bsbX^T \bsbX)}^2  +  K'   J^* \sum_{s=0}^{t-1} \kappa^{t-s}\lambda_s^2.  $$ This allows for     much larger values of  $\lambda_s$ to be used in earlier iterations to attain  the same   accuracy. It              relaxes  the regularity condition required by applying    a fixed threshold level.
\end{remark} 


 At the end, we re-state some results under  $\rho > \| \bsbX\|_2$, to get    more intuition and implications.
For a general  $\bsbX$ (unscaled),     \eqref{recurL1} reads  $$\|  \bsbb^{(t)} - \bsbb^* \|_{(\rho^2 \bsbI - \bsbX^T \bsbX)}^2  \le   \kappa^t\|  \bsbb^{(0)} - \bsbb^* \|_{(\rho^2\bsbI - \bsbX^T \bsbX)}^2  + \frac{\kappa }{1 - \kappa} K'  \sigma^2 \lambda^2 J^*.$$ Set $\rho$  to be a number slightly larger than $\| \bsbX\|_2$, i.e., $\rho=(1+\epsilon) \| \bsbX\|_2$,   $\epsilon>0$. Then, we know that the    prediction error $\|\bsbX \bsbb^{(t)} - \bsbX \bsbb^*\|_2^2$ decays geometrically fast to $O(\sigma^2 J^* \log(e p))$ with high probability, when  $\epsilon$, $\delta$, $\vartheta$, $K$ are viewed as constants;  a similar conclusion is true for the estimation error. This is simply due to  $$ \frac{{\rho^2 - \|\bsbX\|_2^2 }}{{\|\bsbX\|_2^2}}\| \bsbb^{(t)} - \bsbb^*\|_{\bsbX^T \bsbX }^2 \le (\rho^2 - \| \bsbX\|_2^2)   \| \bsbb^{(t)} - \bsbb^*\|_2^{2} \le   \| \bsbb^{(t)} - \bsbb^*\|_{(\rho^2\bsbI - \bsbX^T \bsbX)}^2.$$
Accordingly,     there is no need to run TISP till convergence---one can  terminate the algorithm earlier, at, say,   $t_{\max}= \log\{\rho^2 \|  \bsbb^{(0)} - \bsbb^* \|^2/(K  \sigma^2 \lambda^2 J^* )\}$ $/{\log (1/\kappa)}$, without sacrificing much statistical accuracy. The formula also reflects that the quality of the initial point  affects the required iteration number.

There are some related results in the literature. (i) As mentioned previously, in a broad {convex} setting \cite{agarwal2012} proved the  geometric decay  of the optimization error $\|\bsbb^{(t)} - \hat \bsbb\|$ to the desired statistical precision, where $\hat\bsbb$   is the convergent point.
\cite{Loh15} extended the conclusion to a family of nononvex optimization problems, and they showed  that when some regularity conditions hold, \textit{every} local minimum point   is close to the authentic   $\bsbb^*$. 
In comparison, our results   are derived toward the statistical error between $\bsbb^{(t)}$ and $\bsbb^*$ directly, without requiring   all local minimum points to be statistically accurate. 
(ii) \cite{Zhang10} showed a  similar  fast-converging  statistical error  bound for an elegant multi-stage capped-$\ell_1$
regularization procedure. However,  the procedure carries out an expensive  $\ell_1$ optimization at each step.  Instead,      \eqref{tispdef} involves   a simple and cheap  thresholding, and our analysis covers any $\Theta$.


\section*{Acknowledgement}
The author would like to thank the editor, the associated editor and two anonymous referees
for their careful comments and useful suggestions that
improve the quality of the paper. The author also appreciates Florentina Bunea for the  encouragement.
This work was supported in part by NSF grant DMS-1352259.


\section{Proofs}\label{sec:proofs}
Throughout the proofs, we use $C$, $c$, $L$ to denote universal non-negative constants. They are {not} necessarily the same at each occurrence. Given any matrix $\bsbA$, we use ${\mathcal R}(\bsbA)$ to denote its column space. Denote by    $\Proj_{\bsbA}$  the orthogonal projection matrix onto   ${\mathcal R}(\bsbA)$, i.e., $\Proj_{\bsbA}=\bsbA(\bsbA^T\bsbA)^{+}\bsbA^T$, where $^+$ stands for  the Moore-Penrose pseudoinverse. Let $[p]:=\{1,\cdots,p\}$.
Given $\mathcal J \subset [p]$, we use $\bsbX_{\mathcal J}$ to denote a column submatrix of $\bsbX$ indexed by $\mathcal J$.

\begin{definition}\label{def:subgauss}
$\xi$ is called a sub-Gaussian random variable if there exist constants $C, c>0$ such that $\EP\{|\xi|\geq t\} \leq C e^{-c t^2}, \forall t>0$.  The  scale ($\psi_2$-norm) for  $\xi$ is defined as $\sigma( \xi) =
\inf \{\sigma>0: \EE\exp(\xi^2/\sigma^2) \leq 2\}$. $\bsbxi\in \mathbb R^p$ is called a sub-Gaussian random vector with scale  bounded by $\sigma$ if all one-dimensional marginals $\langle \bsbxi, \bsba \rangle$ are sub-Gaussian satisfying $\|\langle \bsbxi, \bsba \rangle\|_{\psi_2}\leq \sigma \|\bsba \|_2$, $\forall \bsba\in R^{p}$.   \ 
\end{definition}

Examples include Gaussian random variables and bounded random variables such as Bernoulli.
Note that  the assumption that $\vect(\bsbeps)$ is sub-Gaussian does not imply that  the components of $\bsbeps$ must be i.i.d. \\

We begin with two basic facts. Because they are special cases of Lemma 1 and Lemma 2 in \cite{SheGLMTISP}, respectively, we  state them without proofs.
\begin{lemma}
\label{uniqsol-gen-grp}
Given an arbitrary thresholding rule $\Theta$,
let $P$ be any function satisfying $P(\theta;\lambda)-P(0;\lambda)=P_{\Theta}(\theta; \lambda) + q(\theta; \lambda)$ where $P_{\Theta}(\theta; \lambda)\triangleq \int_0^{|\theta|} (\sup\{s:\Theta(s;\lambda)\leq u\} - u) \rd u$, $q(\theta; \lambda)$ is nonnegative and $q(\Theta(t;\lambda))=0$ for all $t$.
Then,  $\hat\beta=\Theta(y;\lambda)$  is always a globally optimal solution to
$
\min_{\beta }  \frac{1}{2}\|y-\beta\|_2^2 + P(|\bsbb|;\lambda).
$
It is the unique optimal solution provided   $\Theta(\cdot;\lambda)$ is continuous at $|y|$.
\end{lemma}

\begin{lemma}
\label{unifuncopt-grp}
Let $Q_0(\beta)= \|y-\beta\|_2^2/2 + P_{\Theta}(|\beta|;\lambda)$. Denote by $\hat\beta$  the unique minimizer of $Q_0(\beta)$. Then for any $\delta$,
$
Q_0(\hat\beta+\delta)-Q_0(\hat\beta) \geq {( 1-{\mathcal L}_{\Theta})} \|\delta\|_2^2/2
$.
\end{lemma}

\subsection{Proof of Theorem \ref{localopt}}
Let $s(u;\lambda):=\Theta^{-1}(u;\lambda) -u$ for $u\ge 0$.
Assume $\hat \bsbb$ is a local minimum point (the proof for a coordinate-wise minimum point follows the same lines). We write $f_{\Theta}$ as $f$ for simplicity. Let $\delta f(\bsbb; \bsbh)$ denote the Gateaux differential of $f$ at $\bsbb$ with increment  $\bsbh$: $\delta f(\bsbb; \bsbh) = \lim_{\epsilon \rightarrow 0+} \frac{f(\bsbb + \epsilon \bsbh)-f(\bsbb )}{\epsilon}$. By the definition of $P_{\Theta}$,  $\delta f(\bsbb, \bsbh)$ exists for any  $ \bsbh \in \mathbb R^p$. Let $l(\bsbb) = \frac{1}{2} \| \bsbX \bsbb - \bsby\|_2^2$. We consider the following directional vectors: $\bsbd_j= [d_1, \cdots, d_p]^T$ with $d_j=\pm 1$ and $d_{j'}=0, \forall j'\ne j$. Then for any $j$,
\begin{align}
\delta l(\bsbb; \bsbd_j) &= d_j \bsbx_j^T (\bsbX \bsbb - \bsby),\\
\delta  P_{\Theta}(\bsbb; \bsbd_j) & = \begin{cases} s(|\beta_j|) \mbox{sgn}(\beta_j) d_j, &\mbox{ if } \beta_j\ne 0, \\ s(|\beta_j|), & \mbox{ if } \beta_j=0.\end{cases}
\end{align}

Due to the local optimality of $\hat\bsbb$,  $\delta f(\hat\bsbb; \bsbd_{j})\ge 0$, $\forall j$.
When $\hat \beta_1 \ne 0$, we obtain $\bsbx_1^T (\bsbX \hat\bsbb - \bsby) + s(|\hat\beta_1|;\lambda) \mbox{sgn}(\hat\beta_1) = 0$. When $\hat \beta_1=0$,  $\bsbx_1^T (\bsbX \hat\bsbb - \bsby) + s(|\hat\beta_1|;\lambda) \ge  0$ and $-\bsbx_1^T (\bsbX \hat\bsbb - \bsby) + s(|\hat\beta_1|;\lambda) \ge  0$, i.e., $|\bsbx_1^T (\bsbX \hat\bsbb - \bsby)| \le  s(|\hat\beta_1|;\lambda)=\Theta^{-1}(0;\lambda)$. To summarize,   when $f$ achieves   a local minimum or a coordinate-wise minimum (or more generally, a \textit{local} coordinate-wise minimum) at $\hat \bsbb$,  we have
\begin{align}
&\hat\beta_j \ne 0 \Rightarrow \Theta^{-1}(|\hat \beta_j|;\lambda)\mbox{sgn} (\hat \beta_j) = \hat \beta_j - \bsbx_j^T (\bsbX \hat\bsbb - \bsby)\label{nzb-1}\\
&\hat\beta_j = 0 \Rightarrow \Theta(\bsbx_j^T (\bsbX \bsbb - \bsby); \lambda) = 0
\end{align}
When $\Theta$ is continuous at $\hat \beta_j - \bsbx_j^T (\bsbX \hat\bsbb - \bsby)$, \eqref{nzb-1} implies that $\hat \beta_j = \Theta(\hat \beta_j - \bsbx_j^T (\bsbX \hat\bsbb - \bsby); \lambda)$. Hence $\hat \bsbb$ must be a $\Theta$-estimator satisfying $\bsbb = \Theta( \bsbb + \bsbX^T \bsby - \bsbX^T \bsbX \bsbb; \lambda)$.

\subsection{Proofs of Theorem \ref{th_oracle} and Theorem \ref{th_oracle-l1type}}
\label{app:subsec:oracleerr}

Given $\Theta$, let $\hat \bsbb$ be any $\Theta$-estimator,  $\bsbb$ be any $p$-dimensional vector (non-random) and $\bsbDelta = \hat \bsbb - \bsbb$.
The first result  constructs a useful criterion for $\hat \bsbb$ on basis of Lemma \ref{uniqsol-gen-grp} and Lemma \ref{unifuncopt-grp}.
\begin{lemma} \label{theta-est-pert}
Any $\Theta$-estimator $\hat \bsbb$ satisfies the following inequality for \textit{\emph{any}} $\bsbb\in \mathbb R^p$
\begin{align}
\begin{split}
&\frac{1}{2} \|\bsbX (\hat \bsbb - \bsbb^*)\|_2^2 + \frac{1}{2} \bsbDelta^T (\bsbX^T \bsbX - \mathcal L_{\Theta}\bsbI) \bsbDelta \\
\le \ & \frac{1}{2} \|\bsbX (\bsbb - \bsbb^*)\|_2^2 + P_{\Theta}(\bsbb; \lambda) - P_{\Theta}(\hat \bsbb; \lambda) + \langle \bsbeps, \bsbX \bsbDelta \rangle,
\end{split}\label{keyineq}
\end{align}
where  $\bsbDelta=\hat \bsbb -  \bsbb$.
\end{lemma}

To handle $\langle \bsbeps, \bsbX \bsbDelta \rangle$, we introduce another lemma.

\begin{lemma} \label{lemma:phostochastic}
Suppose $\| \bsbX\|_2\le 1$ and let $\lambda^o=\sigma\sqrt{ \log (e p)}$. Then there exist universal constants $A_1, C, c>0$ such that for any constants $a\ge 2b>0$, the following event
\begin{align}
\sup_{\bsbb\in \mathbb R^p} \{2\langle \bsbeps, \bsbX  \bsbb \rangle - \frac{1}{ a} \|\bsbX  \bsbb  \|_2^2 - \frac{1}{b} [  P_{H}(\bsbb; \sqrt{ab}{A_1\lambda^o}) ]\} \geq a  \sigma^2 t
\end{align}
occurs with probability at most  $ C \exp(-c t)p^{-c A_1^2}$, where $t\geq 0$.
\end{lemma}

The   lemma plays an important role in   bounding the last stochastic term in \eqref{keyineq}. Its proof is based on the following results.
\begin{lemma} \label{lemma:phcomp}
Suppose $\| \bsbX\|_2\le 1$. There exists a globally optimal solution $\bsbb^o$ to
$
\min_{\bsbb} \frac{1}{2}\| \bsby - \bsbX\bsbb\|_2^2 +  P_{H}(\bsbb; {\lambda})
$
such that for any $j: 1\leq j \leq p$, either $\beta_j^o={0}$ or $|\beta_j^o| \geq \lambda$.
\end{lemma}

\begin{lemma}\label{concenGauss}
Given $\bsbX\in \mathbb R^{n\times p}$ and $J: 1\leq J\leq p$, define $\Gamma_{J} '= \{\bsba\in \mathbb R^{p}: \|\bsba\|_2\leq 1,  \bsba \in {\mathcal R}(\bsbX_{\mathcal J}) \mbox{ for some } \mathcal J: | \mathcal J|=J\}$. Let $P_o'(J) = \sigma^2 \{J+ \log {p\choose J}\}$.
Then for any $t\geq 0$,
\begin{align}
\EP \left(\sup_{\bsba \in \Gamma_{J}'} \langle \bsbeps, \bsba \rangle \geq t \sigma +   \sqrt{L P_o'(J)}\right) \leq C\exp(- ct^2),\label{concbnd1}
\end{align}
where $L, C, c>0$ are universal constants.
\end{lemma}

Let $R= \sup_{1\le J \le p}\sup_{\bsbDelta\in \Gamma_J}\{\langle \bsbeps, \bsbX \bsbDelta \rangle -   \frac{1}{2 b}  P_{H}(\bsbDelta; \sqrt{a b} A_1{\lambda^o}) - \frac{1}{2 a} \|\bsbX  \bsbDelta  \|_2^2\} $, with   $\lambda^o, A_{1}$ given  in Lemma \ref{lemma:phostochastic}. (The starting value of $J$ is   1 because when $J(\bsbDelta) = 0$, $\langle \bsbeps, \bsbX \bsbDelta \rangle=0$.) Substituting it into \eqref{keyineq} gives
\begin{align*}
&\frac{1}{2} \|\bsbX (\hat \bsbb - \bsbb^*)\|_2^2 + \frac{1}{2} \bsbDelta^T (2\bsbX^T \bsbX - \mathcal L_{\Theta}\bsbI) \bsbDelta\\
\le & \frac{1}{2} \|\bsbX (\bsbb - \bsbb^*)\|_2^2 + P_{\Theta}(\bsbb; \lambda) - P_{\Theta}(\hat \bsbb; \lambda) + \frac{1}{2 b}  P_{H}(\bsbDelta; \sqrt{a b} A_1 \lambda^o)  \\& + \frac{1}{2 a} \|\bsbX  \bsbDelta  \|_2^2  + \frac{1}{2} \|\bsbX\bsbDelta\|_2^2+R \\
\le &\frac{1}{2} \|\bsbX (\bsbb - \bsbb^*)\|_2^2 + P_{\Theta}(\bsbb; \lambda) - P_{\Theta}(\hat \bsbb; \lambda) + \frac{1}{2 b}  P_{H}(\bsbDelta; { \sqrt{a b} A_1\lambda^o})\\& + \frac{1}{2 }(1+\frac{1}{a}) \|\bsbX  \bsbDelta  \|_2^2 +R.
\end{align*}
Because $\EP(R\ge a  \sigma^2 t)\le C\exp(-ct)$, we know $\EE [ R] \lesssim a   \sigma^2$.

Let   $\lambda = A\lambda^o$ with $A= A_1 \sqrt{ab}$ and set  $b \ge 1/(2\vartheta)$. The regularity condition ${\mathcal  R}_0(\delta, \vartheta, K, \bsbb, \lambda)$ implies that
\begin{align}    \frac{1}{2 b }  P_{H}(\bsbDelta; {\lambda}) + \frac{\mathcal L_{\Theta}}{2} \|\bsbDelta\|_2^2 \le \frac{2-\delta}{2} \| \bsbX  \bsbDelta\|_2^2 + P_{\Theta}(\hat \bsbb; \lambda)+ K  P_{\Theta}(\bsbb; \lambda). \label{r0inTh} \end{align}
Choose  $a$ to satisfy  $a>1/\delta$, $a\ge 2  b$. Combining the last two inequalities gives
\begin{align}
& \EE [\|\bsbX (\hat \bsbb - \bsbb^*)\|_2^2]\notag  \\ \le & \|\bsbX (\bsbb - \bsbb^*)\|_2^2 +2 (K+1) P_{\Theta}(\bsbb; \lambda)    + \EE[ (1+\frac{1}{a} - \delta) \|\bsbX  \bsbDelta  \|_2^2 ]+2\EE [R] \notag \\ \lesssim  &\|\bsbX (\bsbb - \bsbb^*)\|_2^2 + P_{\Theta}(\bsbb; \lambda) + \sigma^2, \label{r0inTh-2}
\end{align}
with the last inequality due to   $\| \bsbX \bsbDelta\|_2^2 \le (1+1/c) \| \bsbX( \bsbb - \bsbb^*)\|_2^2 + (1+c) \| \bsbX(\hat \bsbb - \bsbb^*)\|_2^2$ for any $c>0$.
\\

The proof of Theorem \ref{th_oracle-l1type} follows the lines of the proof of Theorem \ref{th_oracle}, with  \eqref{r0inTh} replaced by
\begin{align*}
 \frac{1}{2 b }  P_{H}(\bsbDelta; {\lambda}) + \frac{\mathcal L_{\Theta}}{2} \|\bsbDelta\|_2^2 +P_{\Theta}( \bsbb; \lambda)\le \frac{2-\delta}{2} \| \bsbX  \bsbDelta\|_2^2 + P_{\Theta}(\hat \bsbb; \lambda)+ K \lambda^2 J(\bsbb),  \end{align*}
and  \eqref{r0inTh-2} replaced by
\begin{align*}
&\EE [\|\bsbX (\hat \bsbb - \bsbb^*)\|_2^2] \\  \le & \|\bsbX (\bsbb - \bsbb^*)\|_2^2 +2 K \lambda^2 J(\bsbb)  + \EE[ (1+\frac{1}{a} - \delta) \|\bsbX  \bsbDelta  \|_2^2 ]+2\EE [R] \notag\\ \lesssim  &\|\bsbX (\bsbb - \bsbb^*)\|_2^2 + \lambda^2 J(\bsbb) + \sigma^2. 
\end{align*}
The details are omitted.

\subsection{Proof of Theorem \ref{globaloracle}}

From the proof of Lemma \ref{lemma:phcomp}, there exists a $\Theta$-estimator $\hat \bsbb$ which   minimizes  $f(\bsbb)=l(\bsbb) + P_{\Theta}(\bsbb; \lambda)$. This means that    the  term  $\frac{1}{2} \bsbDelta^T (\bsbX^T \bsbX - \mathcal L_{\Theta}\bsbI) \bsbDelta$ can be dropped from  \eqref{keyineq}. Following the lines of  Section \ref{app:subsec:oracleerr},   \eqref{rateoracle} holds under    a modified  version of ${\mathcal  R}_0(\delta, \vartheta, K, \bsbb, \lambda)$, which replaces   \eqref{reg0} with
\begin{align}
    {\vartheta}  P_{H}(\bsbb' - \bsbb; {\lambda}) \le \frac{1-\delta}{2} \| \bsbX  (\bsbb' - \bsbb)\|_2^2 + P_{\Theta}( \bsbb' ; \lambda)+ K  P_{\Theta}(\bsbb; \lambda). \label{reg0-global}
\end{align}
Using the sub-additivity of $P_H$, we know that any design matrix satisfies 
\eqref{reg0-global}    for {any} $0<\vartheta\le 1$,  $\delta\le 1$, $K\ge \vartheta$. 

\subsection{Proof of Theorem \ref{th_seq} and Corollary \ref{corrseq}}
\label{app:subsec:proofth5}

Let  $f(\bsbb) = l(\bsbb) + P_{\Theta}(\bsbb;\lambda)$ where $l(\bsbb) = \frac{1}{2} \| \bsbX \bsbb - \bsby \|_2^2$.

\begin{lemma}\label{lemtri}
Let $\bsbb^{(t+1)} = \Theta( \bsbb^{(t)} + \bsbX^T \bsby - \bsbX^T \bsbX \bsbb^{(t)}; \lambda)$.  Then the following `triangle inequality' holds for any $\bsbb\in \mathbb R^p$
\begin{align*}
& \ \frac{1 - \mathcal L_{\Theta}}{2} \| \bsbb^{(t+1)} - \bsbb \|_{  2}^2 + \frac{1}{2} \| \bsbb^{(t+1)} - \bsbb^{(t)} \|_{\bsbI - \bsbX^T \bsbX}^2   \\
\le & \  \frac{1}{2} \| \bsbb^{(t)} - \bsbb \|_{\bsbI - \bsbX^T \bsbX}^{2}  +f(\bsbb) - f(\bsbb^{(t+1)}).
\end{align*}
\end{lemma}

Letting  $\bsbb = \bsbb^*$ in the lemma, we have
\begin{align*}
&\frac{1}{2} \| \bsbb^{(t+1)} - \bsbb^* \|_{ \bsbX^T \bsbX + (1 - \mathcal L_{\Theta}) \bsbI}^2 + \frac{1}{2} \| \bsbb^{(t+1)} - \bsbb^{(t)} \|_{\bsbI - \bsbX^T \bsbX}^2 + P_{\Theta}(\bsbb^{(t+1)};\lambda) \\
\le &  \frac{1}{2} \| \bsbb^{(t)} - \bsbb^* \|_{\bsbI - \bsbX^T \bsbX}^2 + \langle \bsbeps, \bsbX(\bsbb^{(t+1)} - \bsbb^*)\rangle + P_{\Theta}(\bsbb^{*};\lambda). \end{align*}
Moreover, under ${\mathcal  S}_0(\delta, \vartheta, K, \bsbb^*, \bsbb', \lambda)$ with $\bsbb' = \bsbb^{(t+1)}$,
\begin{align*}
&\vartheta P_H (\bsbb^{(t+1)} - \bsbb^*;\lambda) + \frac{1+\delta} {2} \| \bsbb^{(t+1)} - \bsbb^*\|_2^2 - K P_{\Theta}  (\bsbb^*; \lambda) \\ \le & \frac{1}{2} \|\bsbb^{(t+1)} - \bsbb^* \|_{\bsbX^T \bsbX + (1-\mathcal L_{\Theta})\bsbI}^2  + P_{\Theta} (\bsbb^{(t+1)};\lambda) + \frac{1}{2} \|\bsbb^{(t+1)} - \bsbb^* \|_{\bsbX^T \bsbX}^2.
\end{align*}
Combining the last two inequalities gives
\begin{align*}
&\frac{1+\delta}{2} \| \bsbb^{(t+1)} - \bsbb^* \|_{  \bsbI - \bsbX^T \bsbX }^2 + \frac{1}{2} \| \bsbb^{(t+1)} - \bsbb^{(t)} \|_{\bsbI - \bsbX^T \bsbX} ^2 \\& \quad +\frac{\delta}{2} \| \bsbb^{(t+1)} - \bsbb^* \|_{   \bsbX^T \bsbX }^2 +{\vartheta}  P_{H}(\bsbb^{(t+1)} - \bsbb^{*}; {\lambda})  \\
\le &  \frac{1}{2} \| \bsbb^{(t)} - \bsbb^{*} \|_{\bsbI - \bsbX^T \bsbX}^2+(K+1) P_{\Theta}(\bsbb^{*};\lambda)  + \langle \bsbeps, \bsbX(\bsbb^{(t+1)} - \bsbb^*)\rangle.
\end{align*}

Let  $\Gamma_{J} = \{\bsbb\in \mathbb R^{p}:  J(\bsbb) = J\}$, $\lambda^o=\sigma\sqrt{\log(ep)}$.
 We define an event $\mathcal E$ with its complement given by
$$\mathcal E^c\triangleq\{ \sup_{\bsbb} \{2\langle \bsbeps, \bsbX  \bsbb \rangle - \frac{1}{ a} \|\bsbX  \bsbb  \|_2^2 - \frac{1}{b} [  P_{H}(\bsbb; \sqrt{ab}{A_1\lambda^o}) ]\} \geq 0 \}.$$
By Lemma \ref{lemma:phostochastic},   
there exists a universal constant $L$ such that for any $A_1^2  \ge L$,  $a \ge 2 b>0$,  $P(\mathcal E^{c}) \le C p^{-c A_1^2  }$. Clearly, $\mathcal E$ implies
\begin{align}
\langle \bsbeps, \bsbX(\bsbb^{(t+1)} - \bsbb^*)\rangle \le \frac{1}{2a} \| \bsbb^{(t+1)} - \bsbb^* \|_{   \bsbX^T \bsbX }^2  + \frac{1}{2b} P_H (\bsbb^{(t+1)} - \bsbb^*; \sqrt{ab} A_1 \lambda^o).\label{ineqforcor}
\end{align}

Take $b=1/(2\vartheta)$, $a = 1/(\delta\wedge\vartheta)$, $A_1\ge \sqrt{L}$,  and $\lambda = A_1 \sqrt{a b} \lambda^o$. Then, on $\mathcal E$ we get the desired statistical accuracy bound
 \begin{align*}
 \frac{1+\delta}{2} \| \bsbb^{(t+1)} - \bsbb^* \|_{  \bsbI - \bsbX^T \bsbX }^2 \le\frac{1}{2} \| \bsbb^{(t)} - \bsbb^* \|_{\bsbI - \bsbX^T \bsbX}^2 +(K+1)P_{\Theta}(\bsbb^* ;\lambda).
 \end{align*}

 The bound under  $\mathcal S_1$ can be similarly proved. Noticing that \eqref{ineqforcor}   holds for any $t$,    Corollary \ref{corrseq} is immediately true.

\subsection{Proofs of Lemmas}

\subsubsection{Proof of Lemma \ref{theta-est-pert}} Let $f(\bsbb) = l(\bsbb) + P_{\Theta}(\bsbb;\lambda)$ with $l(\bsbb) = \frac{1}{2} \| \bsbX \bsbb - \bsby\|_2^2$. Define
\begin{align}
g(\bsbb, \bsbg) = l(\bsbb) + \langle \nabla l (\bsbb), \bsbg - \bsbb \rangle + \frac{1}{2} \| \bsbg -\bsbb\|_2^2 + P_{\Theta}(\bsbg;\lambda).\label{gdef}
\end{align}
Given $\bsbb$, $ g(\bsbb, \bsbg)$ can be expressed as
$$
\frac{1}{2} \| \bsbg - (\bsbb - \nabla l(\bsbb))\|_2^2 + P_{\Theta}(\bsbg; \lambda) + c(\bsbb),
$$
where $c(\bsbb)$ depends on $\bsbb$ only.

Let $\hat\bsbb$ be a $\Theta$-estimator satisfying $\hat\bsbb = \Theta(\hat\bsbb - \bsbX^T \bsbX \hat \bsbb + \bsbX^T \bsby; \lambda)$. Based on Lemma \ref{uniqsol-gen-grp} and Lemma \ref{unifuncopt-grp}, we have
$$
g(\hat\bsbb, \hat\bsbb + \bsbDelta) -  g(\hat\bsbb, \hat\bsbb) \ge \frac{1 - \mathcal L_{\Theta}}{2} \| \bsbDelta\|_2^2,
$$
from which it follows that
\begin{align*}
f(\hat\bsbb+\bsbDelta) - f(\hat\bsbb) \ge  \frac{1}{2} \bsbDelta^T (\bsbX^T \bsbX - \mathcal L_{\Theta}\bsbI) \bsbDelta. 
\end{align*}
This holds for any $\bsbDelta\in \mathbb R^p$.
\qed

\subsubsection{Proof of  Lemma \ref{lemma:phostochastic}.}
 Let
\begin{align*}
 l_H(\bsbb) &=  2\langle \bsbeps,  \bsbX  \bsbb \rangle - \frac{1}{ a} \|\bsbX  \bsbb \|_2^2 - \frac{1}{b} [  P_{H}(\bsbb; { \sqrt{a b} A_0  \lambda^o}) ]\\
  l_0(\bsbb) &=  2\langle \bsbeps,  \bsbX  \bsbb \rangle - \frac{1}{ a} \|\bsbX  \bsbb \|_2^2 - \frac{1}{b} [  P_{0}(\bsbb; { \sqrt{a b} A_0  \lambda^o}) ],
\end{align*}
 and   $\mathcal E_H=\{\sup_{\bsbb\in \mathbb R^p } l_H(\bsbb)\geq at\sigma^2\}$, and $\mathcal E_0=\{\sup_{\bsbb\in  \mathbb R^p } l_0(\bsbb)\geq at\sigma^2\}$.
Because $P_0\ge P_H$, $\mathcal E_0 \subset   \mathcal E_H$.
We prove that $\mathcal E_H =   \mathcal E_0$.  The occurrence of  $\mathcal E_H$ implies that $l_H(\bsbb^o) \geq a t \sigma^2$ for any $\bsbb^o$ defined by
\begin{align*}
 \bsbb^o \in \arg \min_{\bsbb}  \frac{1}{ a } \|\bsbX  \bsbb \|_2^2 -2\langle \bsbeps, \bsbX  \bsbb\rangle +\frac{1}{b} [  P_{H}(\bsbb; { \sqrt{a b} A_0  \lambda^o}) ],\end{align*}
 With $a \ge 2b>0$,   Lemma \ref{lemma:phcomp} states that  there exists at least one global minimizer $\bsbb^{oo}$ satisfying $  P_{H}(\bsbb^{oo}; { \sqrt{a b} A_1  \lambda^o})= P_{0}(\bsbb^{oo}; { \sqrt{a b} A_1  \lambda^o})$ and thus $l_H(\bsbb^{oo})=l_0(\bsbb^{oo})$. This means that       $ \sup l_0(\bsbb )\ge l_0(\bsbb^{oo}) =l_H(\bsbb^{oo})\geq a t \sigma^2$. So  $\mathcal E_H \subset \mathcal E_0$, and it suffices to prove $\mathcal E_0^c$ occurs with high probability, or more specifically,  $\EP (\mathcal E_0) \le C \exp(-ct) p^{-c A_1^2}$.

Given   $ 1\leq J \leq p$,   define $\Gamma_{J} = \{\bsbb\in \mathbb R^{p}:  J(\bsbb) = J\}$.
Let $R= \sup_{1\le J \le p} \sup_{\bsbb\in \Gamma_J}\{\langle \bsbeps, \bsbX \bsbb \rangle -   \frac{1}{2 b}  P_{0}(\bsbb; \sqrt{a b} A_1{\lambda^o}) - \frac{1}{2 a} \|\bsbX  \bsbb  \|_2^2\}$.
We will  use  Lemma \ref{concenGauss} to bound its tail probability.

Let $P_o'(J) = \sigma^2 \{J+ \log {p\choose J}\}$. We claim that
\begin{align}
\EP[ \sup_{\bsbb\in \Gamma_J}\{ \langle \bsbeps, \bsbX\bsbb  \rangle - \frac{1}{2a}\|\bsbX \bsbb\|_2^2 - a L P_o'(J) \}> {}a t\sigma^2] \leq C \exp(-c t). \label{firstineqR}
\end{align}
Indeed,
\begin{equation}
  \begin{aligned}
&2 \langle \bsbeps, \bsbX\bsbb\rangle  - \frac{1}{a} \| \bsbX\bsbb\|_2^2 - 2 a L P_o '(J) \\
\le & 2 \langle \bsbeps, \bsbX\bsbb/\|\bsbX\bsbb\|_2\rangle \|\bsbX\bsbb\|_2 - 2 \|\bsbX \bsbb\|_2 \sqrt {L P_o'(J)}  - \frac{1}{2a} \| \bsbX\bsbb\|_2^2 \\
= &  2 \|\bsbX \bsbb\|_2 \left ( \langle \bsbeps, \bsbX\bsbb/\|\bsbX\bsbb\|_2 \rangle -  \sqrt {L P_o'(J)} \right ) - \frac{1}{2a} \| \bsbX\bsbb\|_2^2\\
\le &  2 \|\bsbX \bsbb\|_2 \left( \langle \bsbeps, \bsbX\bsbb/\|\bsbX\bsbb\|_2 \rangle -  \sqrt {L P_o'(J)} \right)_+ - \frac{1}{2a} \| \bsbX\bsbb\|_2^2\\
\le & 2a \left( \langle \bsbeps, \bsbX\bsbb/\|\bsbX\bsbb\|_2\rangle  -  \sqrt {L P_o'(J)} \right)_+^2,
\end{aligned}\label{cauchyin}
\end{equation}
where the last inequality is due to   Cauchy-Schwarz inequality.
 \eqref{firstineqR} now follows from
Lemma \ref{concenGauss}.

Set $A_1 \ge 4 \sqrt L$. We write   $ P_0(\bsbb; \lambda^o)$ with $\bsbb\in \Gamma_J$ as $ P_0(J; \lambda^o)$. Noticing some basic facts that  (i)    $P_o'(J) \le C J \log (ep) \le C P_0(J; \lambda^o)$ due to Stirling's approximation,  (ii)     $ \sqrt {({A_1^2}/{2}) P_0(J; \lambda^o)}  \ge   \sqrt{L P_o'(J)} +    \sqrt{c  A_{1}^2 P_0(J; \lambda^o)}$ for some $c>0$,  and (iii) $J\log (ep) \geq \log p +  J$ for any $J\geq 1$, we get
\begin{align*}
& \EP(R\geq a \sigma^2 t ) \\
\le & \sum_{J=1}^p \EP \left(a  \sup_{\bsbb\in \Gamma_J} \left(\langle \bsbeps, \bsbX\bsbb/\|\bsbX\bsbb\|_2\rangle  -  \sqrt {({A_1^2}/{2}) P_0(J; \lambda^o)} \right)_+^2\ge a  \sigma^2 t\right )\\
= & \sum_{J=1}^p \EP(\sup_{\bsba \in \Gamma_{J}'} \langle \bsbeps, \bsba \rangle-  \sqrt {({A_{1}^2}/{2}) P_0(J; \lambda^o)} \ge   \sigma \sqrt t)\\
 \leq & \sum_{J=1}^p   \EP (\sup_{\bsba \in \Gamma_{J}'} \langle \bsbeps, \bsba \rangle  -   \sqrt{L P_o'(J)} \geq  \sqrt t \sigma + \sqrt{c  A_1^2 P_0(J; \lambda^o)})\\
 \leq & \sum_{J=1}^p  C \exp(-c t) \exp\{- c   A_1^{2} (J+\log (p))\} \\
 \leq & C \exp(-c t) \sum_{J=1}^p \exp(-c A_{1}^2 \log p) \exp(-c  A_{1}^2 J)\\
 \leq & C\exp(-ct ) p^{-c  A_1^2 }   ,
\end{align*}
where   the last inequality due to the sum of  geometric series.
  \qed

\subsubsection{Proof of Lemma \ref{lemma:phcomp}.}
Similar to the proof of Lemma \ref{theta-est-pert}, we set  $f_{H}(\bsbb) = l(\bsbb) + P_H(\bsbb;\lambda)$ with $l(\bsbb) = \frac{1}{2} \| \bsbX \bsbb - \bsby\|_2^2$ and construct $
g_{H}(\bsbb, \bsbg) = f_{H}(\bsbg) + \frac{1}{2} \| \bsbg - \bsbb\|_2^2 -(l(\bsbg) - l(\bsbb) - \langle \nabla l(\bsbb), \bsbg - \bsbb\rangle ).
$
Under  $\| \bsbX\|_2 \le 1$, for any $(\bsbb, \bsbg)$,
$$
g_{H}(\bsbb, \bsbg)  - f_{H}(\bsbg) = \frac{1}{2} (\bsbg - \bsbb)^T (\bsbI - \bsbX^T \bsbX) (\bsbg - \bsbb) \ge 0.
$$

Let $\bsbb^o$ be a globally optimal solution to $\min_{\bsbb}f_H(\bsbb)$. Then $ \bsbg^o:=\Theta_{H}(\bsbb^o - \bsbX^T \bsbX \bsbb^o +\bsbX^T \bsby; \lambda)$ gives
$$
f_{H}(\bsbg^o) \le g_{H}(\bsbb^o, \bsbg^o) \le  g_{H}(\bsbb^o, \bsbb^o) = f_{H}(\bsbb^o),
$$
with the second inequality  due to   Lemma \ref{uniqsol-gen-grp}. Therefore, $\bsbg^o$ must also be a global minimizer of  $f_H$, and by definition, $\bsbg^o$ demonstrates a threshold gap as desired.
\qed

\subsubsection{Proof of Lemma \ref{concenGauss}.}
By definition, $\{\langle \bsbeps,  \bsba\rangle: \bsba \in \Gamma_{J}'\}$ is a stochastic process with sub-Gaussian increments. The induced metric on $\Gamma_{J}'$ is Euclidean: $d(\bsba_1, \bsba_2) = \sigma \|\bsba_1 - \bsba_2\|_2$.

To bound   the metric entropy $\log {\mathcal N}(\varepsilon, \Gamma_{J}', d)$, where ${\mathcal N}(\varepsilon, \Gamma_{J}', d)$ is the smallest cardinality of an $\varepsilon$-net that covers $\Gamma_{J}'$  under  $d$, we notice that $\bsba$ is in a $J$-dimensional ball in $\mathbb R^p$.
The number of such balls $\{\Proj_{\bsbX_{\mathcal J}} \cap B_p(0,1): \mathcal J \subset [p]\}$    is at most  ${p \choose J}$, where $B_p(0,1)$ denotes the unit ball in $\mathbb R^p$.   By a standard volume argument (see, e.g., \cite{VershyninIntro}),
\begin{align}
\log  {\mathcal N}(\varepsilon, \Gamma_{r,J}', d)  \le \log {p \choose J}   (\frac{C \sigma}{\varepsilon})^J=  \log {p \choose J} + J\log ({C \sigma}/{\varepsilon}),
\end{align}
where    $C$ is a universal constant. The  conclusion follows from  Dudley's integral bound \citep{talagrand2005generic}.
\qed

\subsubsection{Proof of Lemma \ref{lemtri}}
We use the notation in the proof of Lemma \ref{theta-est-pert} with $g$ defined in \eqref{gdef}.
By Lemma \ref{uniqsol-gen-grp} and Lemma \ref{unifuncopt-grp}, we obtain $g(\bsbb^{(t)}, \bsbb) - g(\bsbb^{(t)}, \bsbb^{(t+1)})\ge  \frac{1 - \mathcal L_{\Theta}}{2} \| \bsbb^{(t+1)} - \bsbb\|_2^2$, namely,
\begin{align*}
\langle \nabla l (\bsbb^{(t)}), \bsbb - \bsbb^{(t+1)} \rangle + P_{\Theta}(\bsbb) - P_{\Theta}(\bsbb^{(t+1)}) + \frac{1}{2} \| \bsbb - \bsbb^{(t)}\|_2^2 \\ - \frac{1}{2} \| \bsbb^{(t)} - \bsbb^{(t+1)} \|_2^2   \ge  \frac{1 - \mathcal L_{\Theta}}{2} \| \bsbb^{(t+1)} - \bsbb\|_2^2.
\end{align*}
To cancel the first-order term, we give two other inequalities based on  second-order lower/upper bounds: \begin{align*}
&l(\bsbb) - l(\bsbb^{(t)} )- \langle \nabla l(\bsbb^{(t)}), \bsbb - \bsbb^{(t)} \rangle \ge \frac{1}{2} \| \bsbb^{(t)} - \bsbb\|_{\bsbX^T \bsbX}^2,\\
&l(\bsbb^{(t)}) + \langle \nabla l(\bsbb^{(t)}), \bsbb^{(t+1)} - \bsbb^{(t)} \rangle - l(\bsbb^{(t+1)} ) \ge -\frac{1}{2} \| \bsbb^{(t+1)} - \bsbb^{(t)}\|_{\bsbX^T \bsbX}^2.
\end{align*}
Adding the three inequalities together gives    the triangle inequality.
\qed

\bibliographystyle{apalike}
\bibliography{crbib}

\end{document}